\documentclass[a4paper]{article}

\input{glyphtounicode}
\usepackage{mathtools,amssymb}
\usepackage{microtype}

\usepackage{adjustbox}

\usepackage{graphicx}
\graphicspath{{figures/}}
\usepackage{caption,subcaption}

\usepackage{cite}

\usepackage{tabu,booktabs}

\usepackage[pdfusetitle,pdfencoding=auto]{hyperref}

\title{\texorpdfstring{$7x\pm1$: Close Relative of Collatz Problem}{7x\unichar{"00B1}1: Close Relative of Collatz Problem}}
\author{David Barina}
\date{
	\small
	e-mail:~\href{mailto:ibarina@fit.vutbr.cz}{ibarina@fit.vutbr.cz}\\
	\smallskip
	Faculty of Information Technology\\
	Brno University of Technology\\
	Brno, Czech Republic
}

\newcommand\w[1]{\makebox[\ww]{$#1$}}
\def\ww{1.25em}

\begin{document}

\maketitle

\begin{abstract}
We show an iterated function of which iterates oscillate wildly and grow at a dizzying pace.
We conjecture that the orbit of arbitrary positive integer always returns to 1, as in the case of Collatz function.
The conjecture is supported by a heuristic argument and computational results.
\end{abstract}

\noindent
It is conjectured that, for arbitrary positive integer $n$, a sequence defined by repeatedly applying the function
\begin{align}
	f(n) = \begin{cases}
		3n + 1 &:\quad \text{if $n \equiv 1 \pmod{2}$,} \\
		n/2    &:\quad \text{if $n \equiv 0 \pmod{2}$}
	\end{cases}
\end{align}
will always converge to the cycle passing through 1.
The odd terms of such sequence typically rise and fall repeatedly.
The conjecture has never been proven.
The problem is known under several different names, including the Collatz problem, $3x+1$ problem, Syracuse problem, and many others.
There is an extensive literature, \cite{Lagarias2003,Lagarias2006}, on this question.

Its close relative is
\begin{align}
	f(n) = \begin{cases}
		7n + 1 &:\quad \text{if $n \equiv \w{+1} \pmod{4}$,} \\
		7n - 1 &:\quad \text{if $n \equiv \w{-1} \pmod{4}$,} \\
		n/2    &:\quad \text{if $n \equiv \w{ 0} \pmod{2}$,}
	\end{cases}
\end{align}
which also always converges to the cycle passing through 1 when iteratively applied on arbitrary positive integer $n$.
Also here, the odd terms typically rise and fall repeatedly.
It is one of many possible generalizations of the $3x+1$ problem.
However, unlike others, this one shares incredibly many similarities with the original conjecture.

To prove that such sequences always return to 1, one would need to show that these sequences could never repeat the same number twice and they cannot grow indefinitely.
Although the $3x+1$ conjecture has not been proven, there is a heuristic argument, \cite{Tao2011,Lagarias1985,Crandall1978}, that suggests the sequence should decrease over time.
A similar heuristic argument can be used for $7x\pm1$ problem.
The argument is as follows.
If $n$ is odd, then $f(n) = 7n \pm 1$ is divisible by 4; thus two iterations of $f(n) = n/2$ must follow.
Conversely, when $n$ is even, then $f(n) = n/2$ follows.
Furthermore, one can verify that if the input $n$ is uniformly distributed modulo~$2^{l+2}$, then the output of the two branches above is uniformly distributed modulo~$2^l$, for an integer $l \ge 0$.
All branches of the subsequent iteration therefore occur with equal probability.
Now, if the input $n$ is odd, the output of the former branch should be roughly $7/4$ times as large as the input $n$.
Similarly, if the input $n$ is even, the output of the latter branch is $1/2$ times as large as $n$.
If we express the magnitude of $n$ logarithmically, we get expected growth from the input $n$ to the output of the branches above
\[
	\frac{1}{2} \log \frac{7}{4} + \frac{1}{2} \log \frac{1}{2} < 0 \text{.}
\]
Since the growth is negative, the heuristic argument suggests that the magnitude tend to decrease over a long time period.

On positive integers, sequences defined by both the $3x+1$ and the $7x\pm1$ functions eventually enter a repeating cycle $1 \rightarrow \cdots \rightarrow 1$.
When zero is included, there is another cycle $0 \rightarrow 0$ which, however, cannot be entered from outside.
When the $3x+1$ is extended to negative integers, the sequence enters one of a total of three known negative cycles.
These are $-1 \rightarrow \cdots \rightarrow -1$, $-5 \rightarrow \cdots \rightarrow -5$, and $-17 \rightarrow \cdots \rightarrow -17$.
Nevertheless, when the $7x\pm1$ is extended to negative integers, the sequence will always converge to the cycle passing through $-1$.
These cycles are listed in Tables~\ref{tab:3x+1 cycles} and~\ref{tab:7x+1 cycles}.
In contrast to the $3x+1$ problem, every progression in $7x\pm1$ on negative numbers corresponds to negated progression on positive numbers, and vice versa.

\begin{table}[h]
	\centering
	\caption{$3x+1$ problem. Known cycles. Only odd terms due to limited space.}
	\begin{tabu} to \linewidth {X[c]|r}
		\toprule
		\rowfont[c]\bfseries
			{cycle} & {length} \\
		\midrule
			$-17 \rightarrow -25 \rightarrow -37 \rightarrow -55 \rightarrow -41 \rightarrow -61 \rightarrow -91 \rightarrow -17 $ & 18 \\
			$-5 \rightarrow -7 \rightarrow -5$ & 5 \\
			$-1 \rightarrow -1$ & 2 \\
			$+1 \rightarrow +1$ & 3 \\
		\bottomrule
	\end{tabu}
	\label{tab:3x+1 cycles}
\end{table}

\begin{table}[h]
	\centering
	\caption{$7x\pm1$ problem. Known cycles. Only odd terms due to limited space.}
	\begin{tabu} to \linewidth {X[c]|r}
		\toprule
		\rowfont[c]\bfseries
			{cycle} & {length} \\
		\midrule
			$-1 \rightarrow -1$ & 4 \\
			$+1 \rightarrow +1$ & 4 \\
		\bottomrule
	\end{tabu}
	\label{tab:7x+1 cycles}
\end{table}

\begin{figure}
	\begin{subfigure}{\linewidth}
		\centering
		\includegraphics{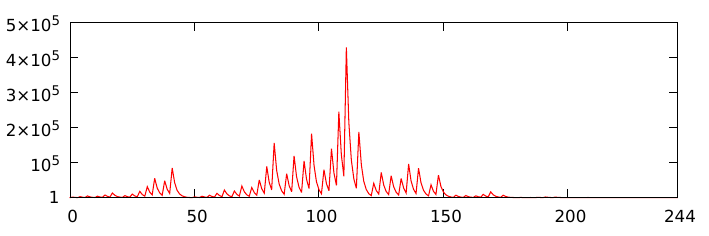}
	\end{subfigure}
	\begin{subfigure}{\linewidth}
		\centering
		\includegraphics{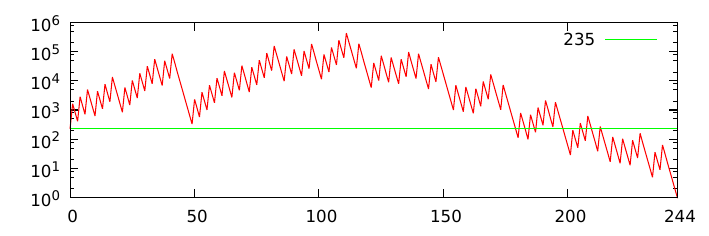}
	\end{subfigure}
	\caption{$7x\pm1$ sequence starting at 235. Due to a very large number range, the sequence in linear scale is shown on top, in logarithmic scale on the bottom.}
	\label{fig:7x+1 235}
\end{figure}

For instance, the $7x\pm1$ sequence for starting value $n=235$ is listed in Table~\ref{tab:7x+1 235}.
It takes 244 steps to reach the number 1 from 235.
This is also known as the total stopping time.
The highest value reached during the progression is 428\,688.
For a better mental picture of this sequence, the progression is also graphed in Figure~\ref{fig:7x+1 235}.
The odd terms can be recognized as local minima, whereas the even terms as either local maxima or descending lines.
One can easily see that the odd terms rise and fall repeatedly.
Such behavior is also common to $3x+1$ sequences.

\begin{table}
	\centering
	\caption{$7x\pm1$ sequence starting at 235. Steps through odd numbers in bold.}
	\begin{adjustbox}{width=\linewidth}
		\small
		\begin{tabu}{p{\textwidth}}
			\toprule
			\textbf{235}, 1644, 822, \textbf{411}, 2876, 1438, \textbf{719}, 5032, 2516,
			1258, \textbf{629}, 4404, 2202, \textbf{1101}, 7708, 3854, \textbf{1927}, 13488,
			6744, 3372, 1686, \textbf{843}, 5900, 2950, \textbf{1475}, 10324, 5162,
			\textbf{2581}, 18068, 9034, \textbf{4517}, 31620, 15810, \textbf{7905}, 55336,
			27668, 13834, \textbf{6917}, 48420, 24210, \textbf{12105}, 84736, 42368, 21184,
			10592, 5296, 2648, 1324, 662, \textbf{331}, 2316, 1158, \textbf{579}, 4052,
			2026, \textbf{1013}, 7092, 3546, \textbf{1773}, 12412, 6206, \textbf{3103},
			21720, 10860, 5430, \textbf{2715}, 19004, 9502, \textbf{4751}, 33256, 16628,
			8314, \textbf{4157}, 29100, 14550, \textbf{7275}, 50924, 25462, \textbf{12731},
			89116, 44558, \textbf{22279}, 155952, 77976, 38988, 19494, \textbf{9747}, 68228,
			34114, \textbf{17057}, 119400, 59700, 29850, \textbf{14925}, 104476, 52238,
			\textbf{26119}, 182832, 91416, 45708, 22854, \textbf{11427}, 79988, 39994,
			\textbf{19997}, 139980, 69990, \textbf{34995}, 244964, 122482, \textbf{61241},
			428688, 214344, 107172, 53586, \textbf{26793}, 187552, 93776, 46888, 23444,
			11722, \textbf{5861}, 41028, 20514, \textbf{10257}, 71800, 35900, 17950,
			\textbf{8975}, 62824, 31412, 15706, \textbf{7853}, 54972, 27486, \textbf{13743},
			96200, 48100, 24050, \textbf{12025}, 84176, 42088, 21044, 10522, \textbf{5261},
			36828, 18414, \textbf{9207}, 64448, 32224, 16112, 8056, 4028, 2014,
			\textbf{1007}, 7048, 3524, 1762, \textbf{881}, 6168, 3084, 1542, \textbf{771},
			5396, 2698, \textbf{1349}, 9444, 4722, \textbf{2361}, 16528, 8264, 4132, 2066,
			\textbf{1033}, 7232, 3616, 1808, 904, 452, 226, \textbf{113}, 792, 396, 198,
			\textbf{99}, 692, 346, \textbf{173}, 1212, 606, \textbf{303}, 2120, 1060, 530,
			\textbf{265}, 1856, 928, 464, 232, 116, 58, \textbf{29}, 204, 102, \textbf{51},
			356, 178, \textbf{89}, 624, 312, 156, 78, \textbf{39}, 272, 136, 68, 34,
			\textbf{17}, 120, 60, 30, \textbf{15}, 104, 52, 26, \textbf{13}, 92, 46,
			\textbf{23}, 160, 80, 40, 20, 10, \textbf{5}, 36, 18, \textbf{9}, 64, 32, 16, 8,
			4, 2, \textbf{1} \\
			\bottomrule
		\end{tabu}
	\end{adjustbox}
	\label{tab:7x+1 235}
\end{table}

\clearpage

\begin{figure}
	\centering
	\begin{subfigure}{\linewidth}
		\centering
		\includegraphics{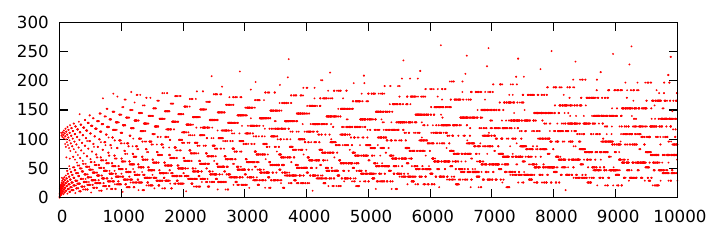}
	\end{subfigure}
	\begin{subfigure}{\linewidth}
		\centering
		\includegraphics{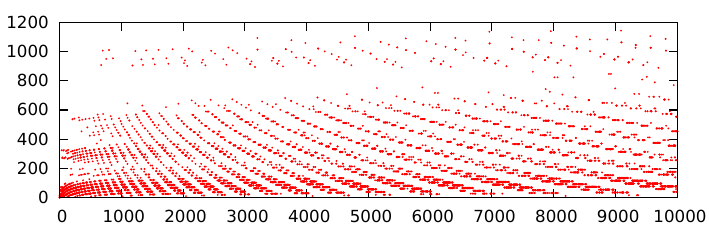}
	\end{subfigure}
	\caption{Numbers 1 to 10\,000 and their total stopping time. The $3x+1$ on top, the $7x\pm1$ on the bottom.}
	\label{fig:total stopping time}
\end{figure}

The progression lengths for both the $3x+1$ and the $7x\pm1$ problems are shown in Figure~\ref{fig:total stopping time}.
Regarding the successive $n$, the behavior of total stopping time is obviously irregular.
Despite this, we can see regular patterns in graphs of these times for both of the problems.
Consecutive starting values tend to reach the same total stopping time.

In order to compare the behavior of the $3x+1$ and $7x\pm1$ sequences, consider following tables.
Tables~\ref{tab:3x+1 longest progression} and \ref{tab:7x+1 longest progression} show the longest progression (total stopping time) for any starting number less than given limit.
One can see that the $7x\pm1$ sequences tend to have recognizably longer progressions.
Moreover, Tables~\ref{tab:3x+1 peak value} and \ref{tab:7x+1 peak value} show that maximum value reached during a progression for any starting number below the given limit.
This value grows significantly faster in the $7x\pm1$ problem that in the $3x+1$ case.

A lot of generalizations, e.g., \cite{Crandall1978,Lagarias1985,Garcia1999,Mignosi1995,Matthews2012}, of the original Collatz problem can be found in the literature.
In \cite{Crandall1978}, the author also mention the $7x+1$ problem.
The definition of such a problem is, however, different from the definition discussed in this paper.
To the best of my knowledge, the $7x\pm1$ function studied in this paper has never appeared before.
I have computationally verified the convergence of the $7x\pm1$ problem for all numbers up to $10^{15}$.

\clearpage

\begin{table}
	\centering
	\caption{$3x+1$ problem. Longest progression for values less than given value.}
	\begin{tabu} to \linewidth {X[c]|X[r]|X[r]}
		\toprule
		\rowfont[c]\bfseries
			{below}   & {peak steps} & {start value} \\
		\midrule
			$10^1$    &     19 &                9 \\
			$10^2$    &    118 &               97 \\
			$10^3$    &    178 &              871 \\
			$10^4$    &    261 &           6\,171 \\
			$10^5$    &    350 &          77\,031 \\
			$10^6$    &    524 &         837\,799 \\
			$10^7$    &    685 &      8\,400\,511 \\
			$10^8$    &    949 &     63\,728\,127 \\
			$10^9$    &    986 &    670\,617\,279 \\
			$10^{10}$ & 1\,132 & 9\,780\,657\,630 \\
		\bottomrule
	\end{tabu}
	\label{tab:3x+1 longest progression}
\end{table}

\begin{table}
	\centering
	\caption{$7x\pm1$ problem. Longest progression for values less than given value.}
	\begin{tabu} to \linewidth {X[c]|X[r]|X[r]}
		\toprule
		\rowfont[c]\bfseries
			{below}   & {peak steps} & {start value} \\
		\midrule
			$10^1$    &     18 &                7 \\
			$10^2$    &    326 &               70 \\
			$10^3$    & 1\,011 &              801 \\
			$10^4$    & 1\,144 &           9\,087 \\
			$10^5$    & 1\,551 &          98\,003 \\
			$10^6$    & 2\,799 &         775\,533 \\
			$10^7$    & 3\,480 &      7\,632\,037 \\
			$10^8$    & 5\,025 &     61\,475\,411 \\
			$10^9$    & 5\,444 &    983\,358\,845 \\
			$10^{10}$ & 5\,717 & 6\,346\,893\,259 \\
		\bottomrule
	\end{tabu}
	\label{tab:7x+1 longest progression}
\end{table}

\clearpage

\begin{table}
	\centering
	\caption{$3x+1$ problem. Maximum value reached in progressions.}
	\begin{tabu} to \linewidth {c|X[4r]|X[r]}
		\toprule
		\rowfont[c]\bfseries
			{below}   &                     {peak value} &    {start value} \\
		\midrule
			$10^1$    &                               52 &                7 \\
			$10^2$    &                           9\,232 &               27 \\
			$10^3$    &                         250\,504 &              703 \\
			$10^4$    &                     27\,114\,424 &           9\,663 \\
			$10^5$    &                 1\,570\,824\,736 &          77\,671 \\
			$10^6$    &                56\,991\,483\,520 &         704\,511 \\
			$10^7$    &           60\,342\,610\,919\,632 &      6\,631\,675 \\
			$10^8$    &       2\,185\,143\,829\,170\,100 &     80\,049\,391 \\
			$10^9$    &  1\,414\,236\,446\,719\,942\,480 &    319\,804\,831 \\
			$10^{10}$ & 18\,144\,594\,937\,356\,598\,024 & 8\,528\,817\,511 \\
		\bottomrule
	\end{tabu}
	\label{tab:3x+1 peak value}
\end{table}

\begin{table}
	\centering
	\caption{$7x\pm1$ problem. Maximum value reached in progressions.}
	\begin{tabu} to \linewidth {c|X[4r]|X[r]}
		\toprule
		\rowfont[c]\bfseries
			{below}   &                                               {peak value} &    {start value} \\
		\midrule
			$10^1$    &                                                         64 &                3 \\
			$10^2$    &                                                   428\,688 &               35 \\
			$10^3$    &                                          20\,492\,891\,264 &              701 \\
			$10^4$    &                                          34\,462\,899\,848 &           8\,317 \\
			$10^5$    &                               965\,557\,666\,410\,854\,560 &          56\,925 \\
			$10^6$    &                           16\,785\,854\,261\,378\,324\,480 &         199\,093 \\
			$10^7$    &                387\,911\,901\,837\,284\,812\,874\,137\,728 &      4\,351\,011 \\
			$10^8$    &           432\,862\,432\,624\,267\,939\,703\,128\,640\,368 &     98\,600\,229 \\
			$10^9$    &   1\,278\,593\,034\,093\,037\,189\,798\,609\,704\,765\,568 &    662\,844\,973 \\
			$10^{10}$ & 421\,614\,662\,439\,923\,712\,249\,655\,593\,962\,998\,304 & 9\,725\,365\,821 \\
		\bottomrule
	\end{tabu}
	\label{tab:7x+1 peak value}
\end{table}

\clearpage

\bibliographystyle{unsrt}
\bibliography{sources}

\end{document}